\documentclass[12pt]{article} 

\usepackage{amsmath,amsfonts,amssymb,amscd}
\let\emptyset\varnothing

\def\rit{\mathbb{R}} 
\def\zit{\mathbb{Z}}   
\def\nit{\mathbb{N}} 
\def\cit{\mathbb{C}} 
\def\fit{\mathbb{F}}

\def\rit{\mathbf{R}} 
\def\zit{\mathbf{Z}}   
\def\nit{\mathbf{N}} 
\def\cit{\mathbf{C}} 

\def\Lc{\mathcal{L}}

\def\Dc{\mathcal{D}}
\def\Hc{\mathcal{H}}
\def\Fc{\mathcal{F}}

\def\Pc{\mathcal{P}}
\def\HH{\mathcal{HH}}
\def\HH{\hbox{{$\mathcal{H}$}\kern-5.2pt{$\mathcal{H}$}}}

\def\Fb{\mathbf{F}}

\def\psib{\mathbf{\psi}}

\def\End{\mathrm{End}}
\def\Res{\mathrm{Res}}

\def\phil{\varphi_{\lambda}}
\def\phim{\varphi_{\mu}}
\def\psil{\psi_{\lambda}}
\def\psim{\psi_{\mu}}

\def\dx{ \frac{\partial}{\partial x} }
\def\st{\star}
\def\la{\langle}
\def\ra{\rangle}

\def\Re{\mathsf{Re}\,}

\def\ftext#1{{\let\thefootnote\relax\footnotetext{\noindent #1}}}

\newcommand{\al}{\alpha}

\newcommand{\qed}{\hfill~~\mbox{$\square$}}

\newtheorem{theorem}{Theorem}[section]
\newtheorem{proposition}[theorem]{Proposition}

\newtheorem{lemma}[theorem]{Lemma}

\begin{document}

\title{\bf Hankel transform via double Hecke algebra} 
\author{ \normalsize Ivan Cherednik $^\dag$, Yavor Markov$^{\star}$}
\date{}
\maketitle  
\begin{center}
{\it
$^{\dag,\,\star,}$ Department of Mathematics, University of
North Carolina,\\
Chapel Hill, NC 27599 -- 3250, USA}\end{center}
\medskip
\centerline{March, 2000}
\medskip

\ftext{ \hspace{-0.6cm}
$^{\dag}$ {\sl E-mail\/{\rm:} chered@math.unc.edu},
Partially supported by NSF grant DMS--9877048\\
\hspace{-0.6cm}$^{\star}$ {\sl E-mail\/{\rm:} markov@math.unc.edu}}

This paper is a part of the course 
delivered by the first author at UNC in 2000.
The focus is on the advantages of the 
operator approach in the theory 
of Bessel functions and the classical Hankel transform.
We start from scratch. 
The Bessel functions were a must for quite a few  
generations of mathematicians
but not anymore. We mainly discuss the {\it master formula} 
expressing the Hankel transform of the product of the 
Bessel function by the Gaussian. 

By the operator approach, we
mean the usage of the Dunkl operator
and the $\HH'',\,$ {\it double H double prime}, 
the rational degeneration of the double affine Hecke algebra.
This includes the transfer from the symmetric theory
to the nonsymmetric one, which is the key tool of the recent development
in the theory of spherical and hypergeometric  functions. 
In the lectures, the Hankel transform was preceded by the standard
Fourier transform, which is of course nonsymmetric, 
and the Harish-Chandra transform, which is entirely symmetric.  

We followed closely the notes of the lectures
not yielding to the temptation of skipping elementary
calculations. We do not discuss the history
and generalizations. Let us give some references.
The master formula is a particular case of that from
[D]. Our proof is mainly borrowed
from [C1] and [C2]. The nonsymmetric  Hankel transform
is due to C. Dunkl (see also [O,J]). We will see that it
is equivalent to the symmetric one, as well as for the  
master formulas (see e.g. [L], 
Chapter 13.4.1, formula (9)).
This is a special feature of the one-dimensional
setup. Generally speaking, there is an implication
nonsymmetric $\Rightarrow$ symmetric, but not otherwise.

We also study the {\it truncated Bessel functions}, which are
necessary to treat negative half-integral $k,$ when the 
eigenvalue of the Gaussian with respect to the
Hankel transform is infinity. They correspond to the
finite-dimensional representations of the double H double
prime, which are completely described in the paper.
We did not find proper references but it is 
unlikely that these functions never appeared before. 
They are very good to demonstrate the operator
technique. 

We thank D. Kazhdan and A. Varchenko, who stimulated the paper
a great deal, and CIME for the kind invitation.

\section{L-operator}
We begin with the classical operator 
$$\Lc=(\dx)^2+\frac{2k}{x}\dx.$$
Upon the conjugation:
\begin{equation}\label{lconj}
\Lc= |x|^{-k}\,\Hc\, 
|x|^{k},\ \
\Hc= (\dx)^2 +\frac{k(1-k)}{x^2}.
\end{equation}
Here $k$ is a complex number.
Both operators are symmetric $=$ even.

The $\varphi$-function is 
introduced as follows:
\begin{equation}\label{l-eigen}
\Lc\phil(x,k)=4\lambda^2\phil(x,k),\qquad
\phil(x,k)=\phil(-x.k),\qquad \phil(0,k)=1.
\end{equation}
We will mainly write $\phil(x)$ instead of $\phil(x,k).$
Since $\Lc$ is a DO of second order, the 
eigenvalue problem has a two-dimensional space of solutions. The even
ones form a one-dimensional subspace
and the normalization
condition fixes $\phil$ uniquely. Indeed,
the operator $\Lc$ preserves the space
of even functions holomorphic at $0.$
The $\phil$ can be of course constructed explicitly, without any
references
to the general theory of ODE.

We look for a solution in the
form $\phil(x,k)=f(x\lambda,k)$. Set $x\lambda=t$. The resulting ODE is
$$ \frac{d^2f}{dt^2}(t) +
2k\frac{1}{t}\frac{df}{dt}(t)-4f(t)=0,\qquad 
\mbox{ a Bessel-type equation. }$$
Its even normalized solution is given by the following series
\begin{equation}\label{feqn}
f(t,k)=\sum_{m=0}^{\infty}\frac{t^{2m}}{m!\,(k+1/2)\cdots (k-1/2+m)}=
\Gamma(k+\frac12)\sum_{m=0}^{\infty}\frac{t^{2m}}{m!\Gamma(k+1/2 +m)}.
\end{equation}
So 
$$f(t,k)=\Gamma(k+\frac12)t^{-k+\frac12}J_{k-\frac12}(2it).$$
The existence and convergence is for all $t\in\cit$ subject to the 
constraint:  
\begin{equation}\label{nhalf}
k\neq -1/2+n,\ n\in \zit_+. 
\end{equation}

The symmetry $\varphi_\lambda(x,k)=\varphi_x(\lambda,k)$ plays a very
important
role in the theory. Here it is immediate. In the multi-dimensional
setup,
it is a theorem.

Let us discuss other (nonsymmetric)
solutions of (\ref{feqn}) and  (\ref{l-eigen}). Looking for 
$f$ in the form $t^\alpha(1+ct+\ldots)$ in a neighborhood
of $t=0,$
we get that the coefficients of the expansion 
$$
f(t)=t^{1-2k}\sum_{m=0}^\infty c_m t^{2m}\hbox{\ at\ } t=0
$$
can be readily calculated from  (\ref{feqn}) and are well-defined
for all $k.$

The convergence is easy to control.
Generally speaking, such $f$
are neither regular nor even. 
To be precise, we get even functions $f$ regular at $0$ 
when $k=-1/2-n$ for an integer $n\ge 0,$ i.e. when
(\ref{nhalf}) does not hold. 
These solutions cannot be normalized as above because
they vanish at $0.$ 

Note that we do not need nonsymmetric $f$ and 
the corresponding $\phil(x)=f(x\lambda)$ in the paper. 
Only even normalized $\varphi$  
will be considered. The nonsymmetric $\psi$-functions
discussed in the next sections
are of different nature.

\begin{lemma}\label{L*-lemma}
 (a) Let $\Lc^\circ$ be the adjoint operator of $\Lc$
with respect to the $\cit$-valued scalar product 
$\la f,g\ra_0=2\int_{0}^{+\infty}f(x)g(x)dx$. 
Then $|x|^{-2k}\Lc^\circ|x|^{2k}=\Lc$.\\
(b) Setting
 $\la f,g\ra=2\int_{0}^{+\infty}f(x)g(x)x^{2k}dx,$ the $\Lc$ is
 self-adjoint with respect to this scalar product, i.e. $\la
 \Lc(f),g\ra=\la f,\Lc(g)\ra$.
\end{lemma}
{\it Proof.} First, the operator multiplication by $x$ is
self-adjoint. Second, $(\dx)^\circ=-\dx$ via integration by parts. 
Finally,
\begin{align}
&x^{-2k}\Lc^\circ 
x^{2k}=x^{-2k}((\frac{d^2}{dx^2})^\circ+(\dx)^\circ(\frac{2k}{x}))x^{2k}=
x^{-2k}((\dx)^2 - \dx(\frac{2k}{x}))x^{2k}\notag\\
&=x^{-2k}(x^{2k}(\dx)^2 +4kx^{2k-1}\dx + 2k(2k-1)x^{2k-2}-
2kx^{2k-1}\dx - 2k(2k-1)x^{2k-2})\notag\\
&=(\dx)^2 +\frac{2k}{x}\dx=\Lc.
\end{align}
Therefore, $\la \Lc(f),g\ra\ =\ 2\int_{0}^\infty\Lc(f)gx^{2k}dx= $
\begin{align}
&2\int_{\rit_+}f\Lc^\circ(x^{2k}g)\,dx
=2\int_{\rit_+}f x^{2k}\Lc x^{-2k}(x^{2k}g)\,dx=\la f,\Lc(g)\ra.
\end{align}

Actually this calculation is not necessary if (\ref{lconj}) is used.
Indeed, $\Hc^\circ=\Hc.$
\qed

\section{Hankel transform}
Let us define the {\it symmetric Hankel transform} on the
space of continuous functions $f$ on $\rit$ such that
$\lim_{x\rightarrow\infty}f(x)e^{cx}= 0$ for any $c\in\rit$.
Provided (\ref{nhalf}),
\begin{equation}
(\fit_k f)(\lambda)=
\frac{2}{\Gamma(k+1/2)}\int_{0}^{+\infty}\phil(x,k)f(x)x^{2k}dx.
\end{equation}
The growth condition makes the transform  well-defined
for all $\lambda\in \cit,$ because 
$$\phil(x,k)\sim \hbox{Const} (e^{2\lambda x}+e^{-2\lambda x})
\hbox{\ at\ }x=\infty.
$$
The latter is standard. 

We switch from $\fit$ on functions to the transform of the
operators: $\fit(A)(\fit(f)=\fit(A(f)).$ Remark that the Hankel
transform of the function is very much different from the transform
of the corresponding multiplication operator.
The key point of the operator technique is the following lemma. 

\begin{lemma} \label{L-hankel}
Using the upper index to denote the variable ($x$ or $\lambda$),
$$(a)\quad \fit(\Lc^x)=4\lambda^2;\qquad
(b)\quad \fit(4x^2)=\Lc^{\lambda};\qquad
(c)\quad \fit(4x\dx)=-4\lambda\frac{d}{d\lambda} - 4 - 8k.$$
\end{lemma}
{\it Proof.} Claim (a) is a direct consequence of
Lemma~\ref{L*-lemma}~(b)
with $g(x)=\phil(x):$
$$\fit(\Lc f)=\la \Lc f,\phil\ra=\la f,\Lc\phil\ra=
4\lambda^2\la f,\phil\ra= 4\lambda^2\fit(f)$$.
Claim (b) results directly from the $x\leftrightarrow \lambda$ symmetry
of $\phi,$ namely, from the relation
$\Lc^{\lambda}\phil(x)=4x^2\psil(x)$. Concerning (c), 
there are no reasons, generally speaking, to expect any simple Fourier transforms
for the operators different from $\Lc.$ However in this particular
case:
$[\Lc^x,x^2]=4x\dx+2+4k$. Appling $\fit$ to both sides and 
using (a), (b),
$[4\lambda^2,\Lc^{\lambda}/4]=\fit(4x\dx) + 2 + 4k$. Finally
$$\fit(4x\dx)=-4\lambda\frac{d}{d\lambda} - 2 - 4k -2 - 4k=
-4\lambda\frac{d}{d\lambda} - 4 - 8k$$.

Note that 
$$[x\dx,x^2]=2x^2,\ [x\dx,\Lc^x]=-2\Lc^x,
$$ 
because operators $x\dx,\Lc$ are homogeneous of degree $2$ and $-2$.
So $e=x^2,$  $f=-\Lc^x/4,$ and $h=x\dx+k+1/2=[e,f]$
generate a representation of the 
Lie algebra $sl_2(\cit).$
\qed

\begin{theorem}({\bf Master Formula}) 
\label{Mth}
Assuming that $\Re k>-\frac12,$
\begin{align}\label{mafo}
 & 2\int_0^{\infty}\phil(x)\phim(x)e^{-x^2}x^{2k}dx=
  \Gamma(k+\frac12)e^{\lambda^2+\mu^2}\phil(\mu),\\
 & 2\int_0^{\infty}\phil(x)\exp(-\frac{\Lc}{4})(f(x))e^{-x^2} x^{2k}dx=
  \Gamma(k+\frac12)e^{\lambda^2}f(\lambda)\notag,
\end{align}
provided the existence of $\exp(-\frac{\Lc}{4})(f(x))$ and the integral
in the second formula.
\end{theorem}   
{\it Proof.} The left-hand side of the first formula
equals $\Gamma(k+1/2)\fit(e^{-x^2}\phim(x))$. We set
$$
\phim^-(x)=e^{-x^2}\phim(x),\ \phim^+(x)=e^{x^2}\phim(x).
$$
They are eigenfunctions of the operators
$$\Lc_-=e^{-x^2}\circ\Lc\circ e^{x^2},\ 
\Lc_+=e^{x^2}\circ\Lc\circ e^{-x^2}.
$$
To be more exact,  $\phim^{\pm}$ is a
unique eigenfunction of $\Lc^{\pm}$ with eigenvalue $2\mu$,
normalized by $\phim^{\pm}(0)=1$.

Express $\Lc_-$ in terms of the operators from the previous lemma.
\begin{align}
e^{-x^2}(\dx)^2e^{x^2}&=e^{-x^2}(e^{x^2}(\dx)^2+2(2x)e^{x^2}\dx+
(2+4x^2)e^{x^2})=(\dx)^2 + 4x\dx + 2+4x^2,\notag\\
e^{-x^2}\frac{2k}{x}\dx e^{x^2}&=e^{-x^2}(e^{x^2}\frac{2k}{x}\dx+
2xe^{x^2}\frac{2k}{x})=\frac{2k}{x}\dx+4k,\notag\\
\Lc_-&=e^{-x^2}((\dx)^2 +\frac{2k}{x}\dx) e^{x^2}=
\Lc+4x\dx+2+4k+4x^2.
\end{align}
Analogously, $\Lc_+=\Lc - 4x\dx-2-4k+4x^2$.
Now we may use Lemma~\ref{L-hankel}:
\begin{align}
\fit(\Lc_-^x)=&\fit(\Lc^x)+\fit(4x^2)+\fit(4x\dx)+\fit(2+4k)=\notag\\
= &4\lambda^2 + \Lc^{\lambda}-4\lambda\frac{d}{d\lambda} - 4 - 8k +2 +
4k
=\Lc_+^{\lambda}.
\end{align}
Thus 
$$L_+^{\lambda}(\fit\phim^-)=\fit(\Lc_-^x)(\fit\phim^-)=\fit(\Lc_-^x\phim^-)=
2\mu\fit(\phim^-),$$ 
i.e. $\fit\phim^- $ is an eigenfunction of $\Lc_+$
with the eigenvalue $2\mu$. Using the uniqueness,
we conclude that $\fit(\phim^-)(\lambda)=$
$C(\mu)e^{\mu^2}\phim^+(\lambda)$
for a constant $C(\mu).$ However the left-hand side of the master
formula
is $\lambda\leftrightarrow \mu$  symmetric as well as 
$e^{\mu^2}\phim^+(\lambda)=
e^{\lambda^2+\mu^2}\varphi_\mu(\lambda).$ 
So  $C(\mu)=C(\lambda)=C.$
Setting $\lambda=0=\mu,$ we get the desired.

The second formula follows from the first for 
$f(x)=\varphi_\mu(x,k).$
Move $\exp(\mu^2)$ to the left to see this. It is linear in terms
of $f(x)$ and  holds for finite linear combinations of $\varphi$ 
and infinite ones provided the convergence. So it is valid for
any reasonable $f.$ We skip the detail.
\qed

\section{Dunkl operator}
The above proof is straightforward. One needs the
self-duality of the Hankel transform and the
the commutator representation for 
$x\,\partial/\partial x.$  The self-duality holds in the general
multi-dimensional theory. The second property is more special. 
Also our proof does not clarify why the master formula is so simple.
There is a ``one-line'' proof of this important formula, which can be
readily generalized. 
It involves the
{ \it Dunkl operator:} 
\begin{equation}
\Dc=\dx -\frac{k}{x}(s-1),\quad \mbox{ where } 
s \mbox{ is the reflection }
s(f(x))=f(-x).
\end{equation}
The operator $\Dc$ is not local anymore, because $s$ is a global
operator
apart from a neighborhood of $x=0.$ We are going to find its 
eigenfunctions. Generally speaking, this may create problems since
we cannot use the uniquness theorems from the theory of ODE.
However everything is  surprisingly smooth.

\begin{lemma}\label{d-prop} Considering  $x$ as the multiplication
operator, 
 \begin{equation}
 s\circ x= -x\circ s, \qquad s\circ \dx=-\dx\circ s,
\end{equation}
$(a)\  \Dc^2=\Lc$ upon the restriction to even functions,\\
$(b)\  s\circ\Dc\circ s=-\Dc$ and $\Dc^2$ fixes the space
 of even functions.
\end{lemma}
{\it Proof.} Indeed, $(s\circ x)(f(x))=s(xf(x))=-xf(-x)=(-x\circ
s)(f(x)).$
The $\dx$ is analogous. Then
\begin{align}\label{dsqr}
  \Dc^2 &= (\dx)^2 -\frac{k}{x}(s-1)\dx - \dx\frac{k}{x}(s-1) + 
  \frac{k}{x}(s-1)\frac{k}{x}(s-1)\notag\\
  & = (\dx)^2 + \frac{k}{x}\dx(s+1) - \dx\frac{k}{x}(s-1) +
 \frac{k}{x}(s-1)\frac{k}{x}(s-1).
\end{align}
It is simple to calculate the final formula but unnecessary. 
Applying (\ref{dsqr}) to symmetric (i.e. even) functions $f(x)$,
the two last terms will vanish, because $(s-1)(f(x))=f(-x)-f(x)=0.$
So $(s+1)(f(x))=f(-x)+f(x)=2f(x),$ and  
$\Dc^2|_{even} = (\dx)^2 + 2\frac{k}{x}\dx=\Lc.$

Claim (b) is obvious. Indeed, $s^2=1,$ $s\dx s=-\dx s^2=-dx$, and 
$ s(\frac{k}{x}(s-1))s=-\frac{k}{x}s(s^2-s)=
-\frac{k}{x}(s^3-s^2)=-\frac{k}{x}(s-1)$.
Thus $s\circ\Dc\circ s=-\Dc$. 

By the way, this implies that
 $s\circ\Dc^2\circ s=\Dc^2$, i.e. $\Dc^2$ commutes with $s.$
So we do not need an explicit formula for $\Dc^2|_{even}$ to see
that it preserves even functions.
\qed

Let us consider the standard  scalar product
$\la f,g\ra_0=\int_{-\infty}^{+\infty}f(x)g(x)dx.$
Here the functions are continuous $\cit$-valued continues on the
real line $\rit.$ One may add the complex conjugation to $g$
but we will not do this. The scalar product is non-degenerate,
so adjoint operators are well-defined. We continue to use
the notation $H^\circ$ for the pairing $\la f,g\ra_0.$  
Let us calculate the adjoint of $\Dc$ with respect to
$|x|^{2k}.$

\begin{proposition}\label{D*-cor} Setting
 $\la f,g\ra=\int_{-\infty}^{+\infty}f(x)g(x)|x|^{2k}dx,$ the Dunkl
 operator  $\Dc$ is  anti self-adjoint with 
 respect to this scalar product, i.e. $\la \Dc(f),g\ra=-\la
 f,\Dc(g)\ra$. Equivalently, $|x|^{-2k}\,\Dc^\circ\,
 |x|^{2k}\,=\,-\Dc$. 
\end{proposition}
{\it Proof.} 
Recall that  $x^\circ=x$ and  $(\dx)^\circ=-\dx,$ where $x$ is
considered as the  multiplication operator. Then $s^\circ=s:$
$$\la s(f),g\ra_0=\int_{-\infty}^{+\infty}f(-x)g(x)dx
=\int_{+\infty}^{-\infty}f(t)g(-t)(-dt)= \la f,s(g)\ra_0
$$ 
for $t=-x.$  Hence,
\begin{align}
|x|^{-2k}\Dc^\circ|x|^{2k}&=|x|^{-2k}(\dx-\frac{k}{x}(s-1))^\circ|x|^{2k}=
|x|^{-2k}((\dx)^\circ - (s-1)^\circ(\frac{k}{x})^\circ)|x|^{2k}\notag\\
&=|x|^{-2k}(-\dx-(s-1)\frac{k}{x})|x|^{2k}=
|x|^{-2k}(-\dx+\frac{k}{x}(1+s))|x|^{2k}\notag\\
 &=|x|^{-2k}|x|^{2k}(-\dx)+|x|^{-2k}(-\frac{2k}{x}|x|^{2k})+
|x|^{-2k}|x|^{2k}\frac{k}{x}(1+s)\notag\\
&=-\dx +\frac{k}{x}(s-1) = -\Dc.
\end{align}
Finally,
\begin{align}
&\la \Dc(f),g\ra = \int_{-\infty}^{+\infty}\Dc(f(x))g(x)|x|^{2k}dx=
\int_{-\infty}^{+\infty}f(x)\Dc^\circ(|x|^{2k}g(x))dx\\
&=\int_{-\infty}^{+\infty}f(x)|x|^{2k}(|x|^{-2k}\Dc^\circ|x|^{2k})(g(x))dx
\notag\\
&=\int_{-\infty}^{+\infty}f(x)(-\Dc(g(x)))|x|^{2k}dx=-\la f,\Dc(g)\ra.
\qquad\qquad\square
\notag
\end{align}

The proposition readily gives that
$|x|^{-2k}\Lc^\circ|x|^{2k}=\Lc$ on even functions $f.$ 
Indeed, 
$$\la \Lc(f),g\ra=\la \Dc^2(f),g\ra = 
\la f,\Dc^2(g)\ra =\la f,\Lc(g)\ra,
$$
provided that $g$ is even too.
Recall that it was not difficult to check this relation
directly. In the multi-dimensional theory,
this calculation is more involved and the usage of
the (generalized) Dunkl operators makes perfect sense.

\section{Nonsymmetric eigenfunctions}
Our next step will be a study of the eigenfunctions of the Dunkl
operator:
\begin{equation}\label{d-eigen}
\Dc\psil(x,k)=2\lambda\psil(x,k),\qquad \psil(0,k)=1.
\end{equation}
We will use the shortcut notation $f^\iota(x)=s(f(x))=f(-x).$

\begin{lemma} \label{D-unique}
There exists a unique solution of the eigenfunction problem
(\ref{d-eigen}) 
for $\lambda\ne 0$. It is
represented in the form
$\psil(x)=g(\lambda x).$
In the case $\lambda=0,$ the solution is given by the formula
$\psi_0=1+Cx|x|^{-2k-1}$, where $C\in\cit$ is an arbitrary constant.
\end{lemma}
{\it Proof.} Assuming that $\psil$ is a solution of (\ref{d-eigen}), let
$$
\psil^0=\frac12(\psil+\psil^\iota),\ \psil^1=\frac12(\psil-\psil^\iota),
$$
be its even and odd parts.
By Lemma~\ref{d-prop}~(b),
$\Dc s (\psil(x))=-s\Dc (\psil(x))=-2\lambda s(\phil(x))$. 
Hence, (\ref{d-eigen}) is equivalent to
\begin{align}\label{d-eigen01}
& \Dc \psil^0=2\lambda\psil^1;\qquad\psil^0(0)=1\notag\\
& \Dc \psil^1=2\lambda\psil^0\qquad\psil^1(0)=0.
\end{align} 
Furthermore,  $\Dc^2\psil^0=4\lambda^2\psil^0$. 
Since $\psil^0$ is even,
$\Lc \psil^0=4\lambda^2\psil$ due to 
Lemma~\ref{d-prop}. Therefore $\psil^0$ has to coincide with
$\phil$ from the first section. This is true for all $\lambda.$
If $\lambda\ne 0,$ 
\begin{equation}\label{psilo}
\psil^1=\frac{1}{2\lambda}\Dc\psil^0=
\frac{1}{2\lambda}\left(\frac{d\psil^0}{dx}-
\frac{k}{x}(s-1)\psil^0\right)=\frac{1}{2\lambda}\frac{d\phil}{dx}.
\end{equation}  
The last equality holds because $\psil^0=\phil$ is even. 
Finally,
\begin{equation}
\psil(x)=\phil(x)+\frac{1}{2\lambda}\phil'(x)=g(\lambda x)
\hbox{\ for\ } g=f+\frac12f',
\end{equation}
where $\phil(x)=f(\lambda x),$ 
$f$ is from (\ref{feqn}),
and $f'$ is the derivative. It is for  $\lambda\ne 0.$

Let us consider the case $\lambda=0$. We have $\Dc\psil^0=0$,
$\psil^0(0)=1$
and $\Dc\psil^1=0$, $\psil^1(0)=0$. Since $\psil^0$ is even, 
$\Dc\psil^0=\frac{d\psil^0}{dx}=0$. Thus $\psil^0(x)=1$. 
The $\psil^1$ is odd. So 
$$\Dc\psil^1(x)=\frac{d\psil^1(x)}{dx}-\frac{k}{x}(s-1)\psil^1(x)=
\frac{d\psil^1(x)}{dx}+\frac{k}{x}(\psil^1(x)-\psil^1(-x))=
\frac{d\psil^1(x)}{dx}+\frac{2k}{x}\psil^1(1). $$
Solving the resulting ODE, $\psil^1(x)=Cx|x|^{-2k-1}$.\qed

In this proof, we used that $\Dc f(x) = f'(x)$ on even functions and 
$$\Dc f(x) = 
f'(x)
+\frac{k}{x}(f(x)-f(-x))=f'(x)+\frac{2k}{x}f(x)=(\dx+\frac{2k}{x})f(x)
$$
on odd functions. By the way,
it makes obvious the coincidence of $\Lc$ with 
$\Dc^2$ on even $f$. Indeed,  $\Dc^2 f(x)=(\dx+\frac{2k}{x})(\dx f(x)).$
For odd $f,$ it is the other way round:
$\Dc^2 f(x)=\Dc(\Dc f(x))=$ $\dx (\Dc f(x))=
\dx(\dx+\frac{2k}{x})(f(x)).$ In particular,
$$\Dc^2 \psil^1(x)=\dx(\dx+\frac{2k}{x})\psil^1(x)=
((\dx)^2+\dx\frac{2k}{x})\psil^1(x).$$
Hence $\psil^1(x)$ is also a solution of a second order differential
equation.
This equation is different from that for $\varphi,$ but not too
different.
Comparing them we come to the important 
definition of the {\it shift operator}.  
We show the dependence of $\Lc$ on $k$ and set 
$\widetilde{\Lc_k}=(\dx)^2+\dx\frac{2k}{x}$. 

\begin{lemma} \label{l-tilde}
(a) $x^{-1}\circ \widetilde{\Lc_k}\circ x = \Lc_{k+1}$.\\
(b) $\widetilde{\Lc_k}\psil^1=4\lambda^2\psil^1$.\\
(c) $(x^{-1}\circ \widetilde{\Lc_k}\circ x)(x^{-1}\psil^1)=
4\lambda^2 (x^{-1}\psil^1)$.
\end{lemma}
{\it Proof.} The first claim:
\begin{align} x^{-1}\circ(\dx)^2\circ x
&=x^{-1}(x(\dx)^2 +2\dx)=(\dx)^2+\frac2x\dx,\notag\\
x^{-1}\circ\dx\frac{2k}{x}\circ x
&=x^{-1}\circ\dx\circ 2k=\frac{2k}{x}\dx,\notag\\
x^{-1}\circ \widetilde{\Lc_k}\circ x& =
(\dx)^2+\frac2x\dx+\frac{2k}{x}\dx=
(\dx)^2+\frac{2(k+1)}{x}dx=\Lc_{k+1}.
\end{align}
Then $\widetilde{\Lc_k}\psil^1
=\Dc^2\psil^1=4\lambda^2\psil^1$ due to
(\ref{d-eigen01}).
Claim (c) is a combination of (a) and (b).\qed

\begin{proposition}({\bf Shift Formula}) \label{D-shift}
\begin{equation}
\frac1x\psil^1(x,k)=\frac{2\lambda}{1+2k}\psil^0(x,k+1), \mbox{ i.e.
}\quad
\frac1x\frac{d\phil}{dx}(x,k)=\frac{4\lambda^2}{1+2k}\phil(x,k+1).
\end{equation}
\end{proposition}
{\it Proof.} Lemma~\ref{l-tilde}~(c) implies that 
$x^{-1}\psil^1(x,k)=C(\lambda,k)\phil(x,k+1)$, because $\phil(x,k+1)$ is
a unique even normalized solution of (\ref{l-eigen}) for $k+1.$
Thanks to (\ref{psilo})
$\psil^1(x,k)=(2\lambda)^{-1}\frac{d\phil}{dx}(x,k)$.
Thus $x^{-1}\phil'(x,k)=C(\lambda,k)\phil(x,k+1)$. 
The constant $C$ readily results from the expantion (\ref{feqn}) of
$\phil(x,k).$ Explicitly:
\begin{align}
0&=(\Lc_k\phil-4\lambda^2\phil)(0)\Rightarrow\notag\\
0& =(2k+1)(x^{-1}\frac{d\phil}{dx})(0,k)-4\lambda^2.
\end{align} 
The shift formula can be of course checked directly
without $\psil^1,$ a good exercise.\qed 

\section{Master formula}
Let us define the {\it nonsymmetric Hankel transform}. 
We consider complex-valued $C^\infty-$ functions $f$ on $\rit$ such that
$\lim_{x\rightarrow\infty}f(x)e^{cx}= 0$ for any $c\in\rit$
and set
\begin{equation}\label{nshankel}
  (\Fc f)(\lambda)=\frac{1}{\Gamma(k+1/2)}
\int_{-\infty}^{+\infty}\psil(x,k)f(x)|x|^{2k}dx
\end{equation}
We assume that  $\Re k>-\frac12$ and always take $\psi_0(x,k)=1.$
Recall that the case $\lambda=0$ is exceptional (Lemma~\ref{D-unique}~):
the dimension of the space of eigenfunctions is $2.$

Let us compute the transforms of our main operators.
Compare it with Lemma~\ref{L-hankel}~: 
it is much more comfortable
to deal with the operators of the first order.
The upper index denotes the variable.
 
\begin{lemma} \label{D-hankel}
$$(a)\quad \Fc(\Dc^x)=-2\lambda;\qquad
(b)\quad \Fc(2x)=\Dc^{\lambda};\qquad
(c)\quad \Fc(s^x)=s^{\lambda}.$$
\end{lemma}
{\it Proof.} The first formula
is an immediate  consequence of Proposition~\ref{D*-cor}~(a)
with $g(x)=\psil(x):$ 
$$\Fc(\Dc f)=\la \Dc f,\psil\ra=-\la f,\Dc\psil\ra=
-2\lambda\la f,\psil\ra= -2\lambda\Fc(f).$$
Claim (b) follows from the $x\leftrightarrow\lambda$ symmetry.
As to (c), use that $\psil(-x)=\psi_{-\lambda}(x).$ 
\qed

\begin{theorem}({\bf Nonsymmetric Master Formula})
\label{NMth}
\begin{align}\label{nsmf}
 & \int_{-\infty}^{\infty}\psil(x)\psim(x)e^{-x^2}|x|^{2k}dx=
  \Gamma(k+\frac12)e^{\lambda^2+\mu^2}\psil(\mu),\\
 &
  \int_{-\infty}^{\infty}\psil(x)\,\exp(-\frac{\Dc^2}{4})(f(x))\,e^{-x^2} 
 |x|^{2k}dx=
  \Gamma(k+\frac12)e^{\lambda^2}f(\lambda).\notag
\end{align}
\end{theorem}   
{\it Proof.} In the first formula, the left-hand side  equals 
$\Gamma(k+1/2)\Fc(e^{-x^2}\psim(x))$. We set
$$\psim^-(x)=e^{-x^2}\psim(x),\ \psim^+(x)=e^{x^2}\psim(x),$$
$$
\Dc_-=e^{-x^2}\circ\Dc\circ e^{x^2},\ 
\Dc_+=e^{x^2}\circ\Dc\circ e^{-x^2}.
$$ The function $\psim^{\pm}$ is an
eigenfunction of $\Dc^{\pm}$ with eigenvalue $2\mu.$
The normalization fixes it uniquely with the standard reservation
about $\mu=0.$

One gets:
$$\Dc_-=e^{-x^2}(\dx - \frac{k}{x}(s-1)) e^{x^2}=\dx +2x -
\frac{k}{x}(s-1)=\Dc + 2x.$$ 
Correspondingly, $\Dc_+=\Dc - 2x$.
Using Lemma~\ref{D-hankel}~,  
$$\Fc(\Dc_-^x)=\Fc(\Dc^x)+\Fc(2x)= -2\lambda + \Dc^{\lambda}=
\Dc_+^{\lambda}.$$ 
Therefore 
$$D_+^{\lambda}(\Fc\psim^-)=\Fc(\Dc_-^x)(\Fc\psim^-)=\Fc(\Dc_-^x\psim^-)=
2\mu\Fc(\psim^-),$$ 
i.e. $\Fc(\psim^-) $ is an eigenfunction of $\Dc_+$
with the eigenvalue $2\mu,$ and $\Fc(\psim^-)(\lambda)=
C(\mu)e^{\mu^2}\psim^+(\lambda)$. Using the 
$\lambda\leftrightarrow\mu$-symmetry,
$C(\mu)=C(\lambda)=C$ and 
$$C=\int_{\rit}e^{-x^2}|x|^{2k}dx=\Gamma(k+\frac12).$$
Cf. the proof of the
symmetric master formula. 

The second formula readily follows from the first provided
the existence of the function
$\exp(-\frac{\Dc^2}{4})(f(x))$ and the corresponding integral.
The latter function has to go to zero at $x=\infty$ faster than
$e^{cx}$ for any $c\in \rit.$ 
\qed

The symmetric master theorem is of course
a particular case of (\ref{nsmf}).
Indeed, we may replace $\psil(x)$ by
$2\phil(x)=\psil(x)+\psil(-x)=\psil(x)+\psi_{-\lambda}(x)$ 
on the left-hand side. Then 
$$\psil(\mu)\mapsto 
\psil(\mu)+\psi_{-\lambda}(\mu)=\psil(\mu)+\psi_{\lambda}(-\mu)=
2\phil(\mu)
$$
on the right-hand side. We use that the factor 
$e^{\lambda^2+\mu^2}$ is even.
Now we either repeat the same transfer for $\mu$ or simply
symmetrize the integrand. 

It is more surprising that the nonsymmetric theorem
can be deduced from the symmetric one. It is a special feature
of the one-dimensional case. Generally speaking,
there is no reason to expect such an implication.
This may clarify why the nonsymmetric Hankel transform 
and $\psi$ were
of little importance in the classical  theory of Bessel
functions. They could be considered as a minor technical improvement
of the symmetric theory. Now we have the opposite point of view.

Let us deduce Theorem~\ref{NMth} from Theorem~\ref{Mth}.
We may assume that  $\lambda,\mu\ne 0.$
Discarding the odd
summands in the integrand, 
\begin{align}
& \int_{\rit}\psil\psim e^{-x^2}|x|^{2k}dx=
\int_{\rit}(\psil^0+\psil^1)(\psim^0+\psim^1) e^{-x^2}|x|^{2k}dx\notag\\
&=\int_{\rit}(\psil^0\psim^0+\psil^1\psim^1) e^{-x^2}|x|^{2k}dx
=\int_{\rit}(\phil\phim+\psil^1\psim^1) e^{-x^2}|x|^{2k}dx.
\end{align}
The integral of $\phil\phim$ 
is nothing else but (\ref{mafo}). Let us use the shift formula to
manage  $\psil^1\psim^1.$ See Proposition~\ref{D-shift}. 
  
We get
$(\psil^1\psim^1)(x,k)=\frac{4\lambda\mu 
x^2}{(1+2k)^2}(\phil\phim)(x,k+1)$ and
\begin{align}
&\int_{\rit}(\psil^1\psim^1)(x,k) e^{-x^2}|x|^{2k}dx=
\frac{4\lambda\mu}{(1+2k)^2}\int_{\rit}(\phil\phim)(x,k+1)
e^{-x^2}|x|^{2(k+1)}dx\notag\\
&=\frac{4\lambda\mu}{(1+2k)^2}\phil(\mu,k+1)e^{\lambda^2+\mu^2}\Gamma(k+\frac32)
=\frac{2\lambda\mu}{(1+2k)}\phil(\mu,k+1)e^{\lambda^2+\mu^2}
\Gamma(k+\frac12)\notag\\
&=\psil^1(\mu,k)e^{\lambda^2+\mu^2}\Gamma(k+\frac12).
\end{align}
This concludes the deduction.
\qed

\section{Double H double prime}
Let $\HH''$ be the double degeneration
of the double affine  Hecke algebra:
\begin{equation}\label{double}
\HH''=\la \,\partial,\, x,\, s\,\mid\, sxs=-x,\quad s\partial s= -
\partial,
\quad [\partial,x]=1+2ks\,\ra .  
\end{equation}

Its {\it  polynomial representation}  $\rho:\HH''\rightarrow \End(\Pc)$
in  $\Pc=\cit[x]$ is as follows:
$$\rho(x)=\mbox{ multiplication by }x,\qquad \rho(s)=s,\qquad
\rho(\partial)=\Dc,$$
where $s$ is the reflection $f\mapsto f^\iota,$ 
$D$ is the Dunkl operator. The first two of the
defining relations of $\HH''$
are satisfied thanks to 
Lemma~\ref{d-prop}. As to the last,
\begin{align}
&[\Dc,x]=
(\dx-\frac{k}{x}(s-1))x-x(\dx-\frac{k}{x}(s-1))\notag\\ 
=&x\dx +1 -x\frac{k}{x}(-s-1)-x\dx+x\frac{k}{x}(s-1)=
1+2ks.
\end{align}

\begin{theorem}\label{faithfull}
(a) Any nonzero finite linear combination
$H=\displaystyle\sum_{m,n,\epsilon}c_{m,n,\epsilon}x^m\partial^n
s^{\epsilon}$, where $n,m\in\zit_+$, $\epsilon\in\{0,1\}$, acts as a
nonzero operator in $\Pc$.\\
(b) Any element $H\in\HH''$ can be uniquely expressed in the form
$H=\displaystyle\sum_{m,n,\epsilon}x^m\partial^n
s^{\epsilon}.$ The representation $\rho$ is faithful for any $k.$
\end{theorem}
{\it Proof.} 
Let $H=\displaystyle{\sum_{n,\epsilon}f_n\partial^n
s^{\epsilon}=\sum_{m=0}^N(g_{m}\partial^m)(1+s)+
\sum_{m=0}^M(h_{m}\partial^m)(1-s)}$ 
for polynomials $f_m,g_m,$ and $h_m.$ We may assume that
and at least one of the leading coefficients
$g_N(x),\ $ $h_M(x)$ is nonzero.
Then $\rho(H)=L_+(1+s)+L_-(1-s)$ for  differential operators 
$L_+=g_{N}(x)(\dx)^N+\ldots$ and 
$L_-=h_M(x)(\dx)^M+\ldots$ modulo
differential operators of lower orders.
Applying $\rho(H)$ to even and odd functions, we get that 
$\rho(H)=0$ implies that both $L_+$ and $L_-$ have infinite
dimensional spaces of eigenfunctions. This is impossible.
Claim (a) is verified.

Concerning (b),
any element  $H\in \HH''$ can be obviously expressed in the desired
form. Such expression  is unique and the representation
$\rho$ is faithful thanks to (a).\qed  

The theorem is the key point of the representation theory of
the double H. It is a variant of the so-called PBW theorem.
There are not many algebras in mathematics and physics possessing
this property. All have  
important applications. The double Hecke algebra is one of them.  

Next, we study the irreducibility of $\rho$.
\begin{lemma}\label{d-poly-eigen}
The Dunkl operator $\Dc$ has only one eigenvalue in $\Pc$,
namely, $\lambda=0$.
If $k\ne -1/2-n$ for any $n\in\zit_+$, then $\Dc$ has a unique (up
to a constant) eigenfunction in $\Pc$, the constant function $1$. 
When $k= -1/2-n$ for  $n\in\zit_+$,  the space of
$0$-eigenfunctions is $\cit+\cit x^{2n+1}$.
\end{lemma}
{Proof.} Let $p(x)\in \Pc$ be an eigenfunction for $\Dc$.
Since $\Dc$ lowers the degree of any polynomial by 1, we have
$\Dc^{m+1}p=0$, where $m=\deg p$. Therefore all eigenvalues of $\Dc$ are
zero. Representing $p$ as the sum
$p(x)=p^0(x)+p^1(x)$ of even $p^0$ and odd $p^1,$
$\Dc p=\dx p^0+(\dx p^1 +\frac{2k}{x}p^1)=0$. Both
the even and the odd parts of this expression have to be zero. Hence
$\dx p^0=0$ and $\dx p^1+\frac{2k}{x}p^1=0$. Therefore $p^0=$Const.
Setting $p^1(x)=\sum a_l x^{2l+1},$
$$(\dx+\frac{2k}{x})p^1(x)=
\sum a_l(2l+1+2k)x^{2l}=0.$$ 
If $k\ne -1/2-n$ for any $n\in\zit_+$
then $a_l=0$ for any $l$, i.e $p^1=0.$ 
Otherwise $k= -1/2-n$ for a certain
$n\in\zit_+$ and $p^1(x)$ is proportional to $x^{2n+1}$.\qed

\begin{theorem}\label{irrep}\hfill\\ 
(a) The representation $\rho$ is irreducible if and only
if $k\ne -1/2-n$ for any $n\in\zit_+$.\\
(b) If $k= -1/2-n$ for  $n\in\zit_+$, then there exists a unique
non-trivial $\HH''\,$-submodule of $\Pc$, namely, 
$(x^{2n+1})=x^{2n+1}\Pc$.
\end{theorem}
{\it Proof.} Let $\{0\}\ne V\subset\Pc$ be a 
$\HH''$ submodule of $\Pc$. 
Let $0\neq v\in V$ and $m=\deg v$. Set $\Pc^{(m)}=\{p\in\Pc\,|\,
\deg p\leq m\}$. We have $V^{(m)}=V\cap\Pc^{(m)}\ne \{0\}$. 
Then $V^{(m)}$ is $\Dc$ invariant. More exactly,
$\Dc(V^{(m)})\subset V^{(m-1)}.$ 
Thus it
contains an eigenfunction $v_0$ of $\Dc$. If $k\ne -1/2-n$ for  
$n\in\zit_+$ then Lemma~\ref{d-poly-eigen} implies that $v_0=1$. So
$1\in V$ and $V\ni\rho(x^m)(1)=x^m$ for any
$m\in\nit.$ This means that $V=\Pc$ and  completes  (a).

If $k= -1/2-n$ for  $n\in\zit_+$ then Lemma~\ref{d-poly-eigen} states
that $v_0=c_1+c_2 x^{2n+1}\in V$ for constants $c_1,c_2.$
If $c_1\neq 0$ then $s(v_0)+v_0=2c_1\in V$ (the latter is $s$-invariant
and  $V=\Pc$ as above. If $v_0=x^{2n+1}\in V$ then this argument gives
that $(x^{2n+1})=x^{2n+1}\Pc\subset V.$ Moreover $V$ cannot contain
the polynomials of degree less than $2n+1.$ Otherwise we can find
a $0$-eigenvector of $\Dc$ in the space of such polynomials, which
is impossible. Hence $V=(x^{2n+1}).$

The latter is invariant with respect to $x$ and $s.$
Its $\Dc$-invariance readily follows from the formula 
$\Dc(x^l)=(l+(1-(-1)^l)k)x^{l-1}$ considered
in the range  $l\ge 2n+1.$
\qed

We can reformulate the theorem as follows. The 
polynomial representation has a nontrivial (proper) $\HH''\,$-quotient
if
and only if  $k= -1/2-n$ for  $n\in\zit_+.$ In the latter case, such
quotient
is unique, namely, $V_{2n+1}=\Pc/(x^{2n+1}).$ Its dimension is
$2n+1.$ 

Note that the subspace $V^0_{2n+1}$ of $V_{2n+1}$ 
generated by even
polynomials is invariant with respect to the action
of $h=x\dx+k+1/2$ and $e=x^2,$ 
$f=-\Lc/4,$ satisfying the defining relations
of $sl_2(\cit)$ (see Section 2). 
We get an irreducible representation of $sl_2(\cit)$
of dimension $n+1.$

\section{Algebraization}
Let us use $\HH''$ to formalize the previous considerations and to
switch to the standard terminology of the representation theory.

{\bf  (a) Inner product.}  We 
call a representation $V$ of $\HH''$ {\it pseudo-unitary} if it
possesses
a non-degenerate $\cit$-bilinear form $(u,w)$ such that
$(Hu,w)=(u,H^\st v)$ for $H\in \HH''$ for the anti-involution
\begin{equation}\label{istar}
  \partial^{\st}=-\partial,\qquad s^{\st}=s,\qquad x^{\st}=x.
\end{equation}
By anti-involution, we mean a $\cit$-automorphism
satisfying $(AB)^\st=B^\st A^\st.$ We call such form 
{\it $\st$-invariant}.
Formulas (\ref{istar}) are compatible with the defining relations
(\ref{double})
of $\HH''$ and therefore can be extended to the whole $\HH''.$ 
This is straightforward. For instance,
$$ [x^{\st},\partial^{\st}]=[x,-\partial]=[\partial,x]=1+2ks=
1+2ks^{\st}=[\partial,x]^{\st}.$$
We add ``pseudo'' because the pairing, generally speaking,
is not supposed to be positive and the functions can be complex-valued. 

The pairing $\la f,g\ra=\int_{\rit}f(x)g(x)|x|^{2k}dx$ gives an example, 
provided the existence of the integral. 
Taking real-valued functions, we make this inner product 
positive (no ``pseudo'').
Assuming that the functions
are $C^\infty,$ we need to examine the convergence at $x=0$ and
$x=\infty.$  
If $\Re(k)>-1/2$ then it suffices to take 
regular $f$ at $x=0.$ At infinity, $f(x)|x|^{k}$ has to be
of type $L^1(\rit).$ Polynomials times the Gaussian $e^{-x^2}$ are fine.

{\bf  (b) Gaussians.}  A homomorphism $\gamma:V\to W$ for two
$\HH''$-modules
$V,W$ is called a {\it Gaussian} 
if $\gamma H =\tau(H)\gamma$ for the following automorphism  $\tau$ of
$\HH'':$
\begin{equation}
\tau(\partial)=\partial-2x,\qquad \tau(x)=x,\qquad \tau(s)=s.
\end{equation} 
These formulas can be extended to an automorphism of $\HH''.$
Indeed, 
$\tau(s)\tau(\partial)\tau(s)=s(\partial+2x)s=-\partial-2x=-\tau(\partial),$
the same holds for $x$, and
$$
[\tau(\partial),\tau(x)]\ =\ [d-2x,x]\ =\ 1+2s\ =\ \tau(1+2s).
$$

Note that $2$ can be replaced by any constant $\al\in \cit$ in this
definition.
We get a family of automorphisms $\tau_\al(\partial)=
\partial-2\al x$ of
$\HH''.$ 
They lead to the
following  generalization of the master formula: 
$$\int_{\rit}\psil(x)\psim(x)e^{-\alpha x^2}|x|^{2k}dx=
\psil(\frac{\mu}{\alpha})
e^{\frac{\lambda^2+\mu^2}{\alpha}}\alpha^{-k}
\Gamma(k+\frac12).$$
Here $\al>0$ to ensure the convergence. The substitution $u=\sqrt{\alpha}x$
readily makes it equvalent to (\ref{nsmf}): use that 
$\psi_\lambda(x)$
is a function of the $x\lambda.$ One can also follow the proof of the master
formula employing 
$e^{-\alpha x^2}\Dc e^{\alpha x^2}=\Dc+2\alpha x.$ 

If representations $V,W$ are algebras of functions on the same set then 
$\gamma$ can be assumed to be a function, to be more exact, the operator
of 
multuplication by a function. For instance, the multiplication by
$e^{x^2}$ 
sends the polynomial representation $\Pc$ to $\Pc e^{\pm x^2}.$
The latter is a $\HH''$-module too.   
Adding all integral powers of $e^{x^2}$ to $\Pc$ we make this
multiplication
an inner automorphism of the resulting algebra.
However it is somewhat artificial. Algebraically, 
the resulting representation $\Pc[ e^{mx^2},\,m\in\zit]$ is ``too''
reducible.
Analytically, we mix together $e^{x^2}$ and $e^{-x^2},$ functions
with absolutely different behaviour at infinity. There are interesting
examples of
inner automorphisms $\tau,$ but they are finite-dimensional.

{\bf  (c) Hankel transform.}  Following Lemma~\ref{D-hankel}, the 
{\it operator
Hankel transform} is the following automorphism of $\HH'':$
\begin{equation}\label{ohankel}
 \sigma(s)=s,\qquad \sigma(\partial)=-2x,\qquad \sigma(2x)=\partial.
\end{equation}
These realations are obviously comapatible with the
defining relations of $\HH''.$ Any homomorphism
$\Fc: V\to W$ of $\HH''$-modules inducing $\sigma$ on $\HH''$ can
be called a Hankel transform. The main example so far is 
$\Fc:\Pc e^{-x^2}\to \Pc e^{+x^2},$ 
where we identify $x$ and $\lambda$ in
\ref{nshankel}. Indeed,
$\Fc\circ\Dc=-2x\circ \Fc$, $\Fc\circ 2x=\Dc\circ\Fc$, 
$\Fc\circ s=s\circ \Fc$ upon this identification.

It is interesting to interpret the master formula from this viewpoint.
It is nothing else but  the following
identities for $\Fb=  e^{x^2} e^{\partial^2/4}e^{x^2}$:
\begin{equation}\label{omaster}
 \Fb s=s\Fb,\qquad \Fb\partial=-2x\Fb,\qquad \Fb(2x)=\partial\Fb.
\end{equation}
This means that  $\Fb$ is Hankel transform whenever it is well-defined.
Relations (\ref{omaster}) can be deduced directly from the defining
relations.
In the first place, check that $[\partial,x^2]=2x,
[\partial^2,x]=2\partial.$  Then
get that  $[\partial,e^{x^2}]=2xe^{x^2},$ $[e^{\partial^2/4},x]=$
$\partial e^{d^2}/2$ and use it
as follows:
\begin{align}
&e^{x^2} e^{\partial^2/4}e^{x^2} 2x= 
e^{x^2} e^{\partial^2/4}(2x) e^{x^2}=
e^{x^2}(\partial+2x) e^{\partial^2}e^{x^2}=\notag\\
&(\partial -2x+2x)e^{x^2}e^{d^2}e^{x^2}=
\partial e^{x^2} e^{\partial^2}e^{x^2},\notag\\
&e^{-x^2} e^{-\partial^2/4}e^{-x^2} 2x= 
e^{-x^2} e^{-\partial^2/4}(2x) e^{-x^2}=
e^{-x^2}(-\partial+2x) e^{-\partial^2}e^{-x^2}=\notag\\
&(-\partial-2x+2x)e^{x^2}e^{d^2})e^{x^2}=
-\partial e^{x^2} e^{\partial^2}e^{x^2}.\notag
\end{align}
The commutativity of $\Fb$ with $s$ is obvious.

Note the following {\it braid identity} which can be proved
similarly:
$$
\Fb=  e^{x^2} e^{\partial^2/4}e^{x^2}=
=  e^{\partial^2} e^{x^2/4}e^{\partial^2}.
$$

Actually we do need calculations from scratch because
it suffices to use the nonsymmetric master formula
Lemma~\ref{D-hankel} and the fact that
$\Pc$ is a faithful representation. For instance,
$\Fb(2x)\Fb^{-1}$ and $\partial$ coincide in
$\Pc.$ The previous consideration shows that
the former is an element of $\HH''.$ Hence they must
coincide identically, i.e. in $\HH''.$

A good demonstration of the convenience of such an 
algebraization will be the case of negative half-integral
$k.$ Before switching to this case, let us conclude the
``analytic'' theory  
calculating the inverse Hankel transform.

\section{Inverse transform and Plancherel formula}

Let $\Re k>-1/2.$ We use $\psil(x)$ from (\ref{d-eigen}).     
\begin{align}
(\Fc_{re} f)(\lambda)&=\frac{1}{\Gamma(k+1/2)}
\int_{-\infty}^{+\infty}\psil(x)f(x)|x|^{2k}dx,\\
(\Fc_{im} g)(x)&=\frac{1}{\Gamma(k+1/2)}
\int_{-i\infty}^{+i\infty}\psi_x(-\lambda)g(\lambda)
|\lambda|^{2k}d\lambda.
\end{align}  
The first is nothing else but $\Fc$ from (\ref{nshankel}).
We just show explicitely that the
integration is real.

Here we may consider $\cit$-valued functions $f$ on $\rit$ and
$g$ on $\cit$ respectively such that
$$
f(x) = o(e^{cx}) \mbox{ at } x=\infty 
\, \forall c\in \rit\hbox{\ \,and\ \,} f\in
g(\lambda) = o(e^{ci\lambda}) \mbox{ at } \lambda=i\infty
\, \forall c\in \rit.
$$ 
Restricting ourselves with the polynomials times the Gaussian,
$$
\Fc_{re}:\cit[x]e^{-x^2}\rightarrow \cit[\lambda]e^{\lambda^2},\ 
\Fc_{im}:\cit[\lambda]e^{\lambda^2}\rightarrow \cit[x]e^{-x^2}.
$$
The first map is an isomorphism. Let us discuss the latter.

Let $p(x)\in\cit[x]$. Applying the master formula to
$p_i(x)=p(ix),$
$$\frac{1}{\Gamma(k+1/2)}\int_{-\infty}^{+\infty}
\psil(x)p(ix)e^{-x^2}|x|^{2k}dx=\Fc_{re} (e^{-x^2}p(ix))=
e^{\lambda^2}\exp((\Dc^{\lambda})^2/4)(p(i\lambda)).
$$
Since $\Dc\circ I=iI\circ \Dc$ for $I(f)(x)=f(ix),$
$$
\frac{(\Dc^{\lambda})^2}{4} p(i\lambda)=
\left(\frac{-(\Dc^u)^2}{4} p(u)\right)|_{u=i\lambda}.
$$ 

Now we replace $\lambda\mapsto i\lambda,$ use that
$\psi_{i\lambda}(x)=\psil(ix),$ and then integrate by 
substitution using $z=ix.$
The resulting formula reads:
\begin{align}
&\frac{1}{\Gamma(k+1/2)}\int_{-\infty}^{+\infty}
\psil(ix)p(ix)e^{-x^2}|x|^{2k}dx=
e^{-\lambda^2}\exp(-(\Dc^\lambda)^2/4)(p(-\lambda))\\
&=\frac{1}{i\,\Gamma(k+1/2)}\int_{-i\infty}^{+i\infty}        
\psil(z)p(z)e^{z^2}|z|^{2k}dz.
\end{align}
Switching to $\lambda,$
$\Fc_{im}(e^{\lambda^2}p(\lambda))=
e^{-x^2}\exp(-(\Dc^x)^2/4)(p(x)).$

We come to the inversion theorem.

\begin{theorem}({\bf Inversion Formula})
$\Fc_{im}\circ\Fc_{re}=\mathbf{id}\ $ in 
$\cit[x]e^{-x^2},\ $
$\Fc_{re}\circ\Fc_{im}=\mathbf{id}\ $ in 
$\cit[\lambda]e^{\lambda^2}.$
\end{theorem}
{\it Proof:}
$\Fc_{im}\circ\Fc_{re} (e^{-x^2}p(x))
=e^{-x^2}\exp(-\Dc^2/4)\exp(\Dc^2/4)(p(x))=e^{-x^2}p(x).$
The second formula is analogous.
\qed

There is a simple
``algebraic'' proof based on the facts that the 
transform $\Fc_{im}\circ\Fc_{re}$ sends $\Dc\mapsto\Dc$, $2x\mapsto 2x$,
$s\mapsto s.$ Thanks to the irreducibility of  
$\cit[x]e^{\pm x^2},$ we may apply the Schur lemma.
However the spaces are infinite
dimensional, so a minor additional consideration is necessary.
We will skip it because it practically coincides with that
from the proof of the Plancherel formula. 
 
Provided $\Re k>-\frac12,$ the inner products
\begin{equation}
 \la f,g \ra_{re}=\int_{-\infty}^{+\infty}f(x)g(x)|x|^{2k}dx,\qquad
\la f,g \ra_{im}=\frac1i\int_{-i\infty}^{+i\infty}f(\lambda)
g(-\lambda)|\lambda|^{2k}dx
\end{equation}
are non-degenerate respectively in $\cit[x]e^{-x^2}$
and $\cit[\lambda]e^{\lambda^2}$. 

It is obvious when $\rit[x],\rit[\lambda]$
are considered instead of $\cit[x],\cit[\lambda]$ and $k\in \rit.$
Indeed, both forms become positive in this case.
Concerning the second, use that $
g(-\lambda)=\overline{g(\lambda)}$ for a real polynomial $g.$

For complex-valued functions and $k\in \cit,$ the claim requires
proving. Let us use the irreducibility of $\HH''-$ modules 
$\cit[x]e^{-x^2}$ and $\cit[x]e^{x^2},$ which is equivalent
to the irreducibility of the polynomial representation $\Pc=
\cit[x],$ which we already know for $\Re k>-\frac12.$
Then the radical of the form
$\la\,.\,,\,.\,\ra_{re}$ is a submodule of $\cit[x]e^{-x^2}.$ 
It does not coincide with the whole space, 
since 
\begin{equation}\label{gaure}
\la e^{-x^2},e^{-x^2}\ra_{re}
=\int_{\rit}e^{-2x^2}|x|^{2k}dx =(\sqrt2)^{-2k-1}\Gamma(k+\frac12).
\end{equation}
The same argument works for the imaginary integration.

\begin{theorem} ({\bf Plancherel Formula})
\begin{equation}\label{planch} 
\la f,g \ra_{re}=\la \widehat{f},\widehat{g}\ra_{im}
\mbox {  for all  } f,g\in \cit[x]e^{-x^2},\ 
\widehat{f}=\Fc_{re}(f),\,\widehat{g}=\Fc_{re}(g).
\end{equation}
\end{theorem}
{\it Proof.} Setting $\Pc_-=\cit[x]e^{-x^2},$ we need to check that
$\la\,.\,,\,.\,\ra=\la\,.\,,\,.\,\ra_{re}$ coincides with
$\la f,g \ra_1=\la\widehat{f},\widehat{g}\ra_{im}$ for all
$f,g\in\Pc_-$. In the first place,
$\la Hf,g\ra_1=\la f,H^{\st}g\ra_1$ for any
$H\in\HH'',$ i.e. this bilinear form is $\st$-invariant.
Indeed,
\begin{align}
&\la \Dc f,g\ra_1=\la\widehat{\Dc f},\widehat{g}\ra_{im}=
\la -2x\widehat{f},\widehat{g}\ra_{im}=
-\la\widehat{f},-2(-x)\widehat{g}\ra_{im}=
-\la\widehat{f},\widehat{\Dc g}\ra_{im}=-\la f,\Dc g\ra_1.\notag
\end{align}
Similarly, $\la 2x f,g\ra_1=\la f,2xg\ra_1$ and
$\la s(f),g\ra_1=\la f,s(g)\ra_1.$ 

Setting $\la\,.\,,\,.\,\ra_2=\la\,.\,,\,.\,\ra-
\la\,.\,,\,.\,\ra_1,$ we get $\la e^{-x^2},e^{-x^2}\ra_2=0.$ 
Cf. \ref{gaure}. It is a $\star$-invariant
form as well. Let us demonstrate that it vanishes identically.

First, $\la e^{-x^2},(\Dc-2x)f\ra=
\la -(\Dc+2x)e^{-x^2},f\ra=\la 0,f\ra=0$ for any $f \in \Pc_-$
due to the $\st$-invariance. So it is applicable to
$\la\,.\,,\,.\,\ra_2$ too. Second, 
$\cit e^{-x^2}\cap (\Dc-2x)\Pc_-=\emptyset$ because
$\la e^{-x^2},e^{-x^2}\ra\neq 0.$ Third, 
$\Pc_-=\cit e^{-x^2}\oplus (\Dc-2x)\Pc_-.$ Really,
the dimension of  $(\Dc-2x)\Pc_{n-}$ for $ \Pc_n=$
$\cit[1,x,\ldots,x^n]e^{-x^2}$ is $n+1$ since the kernel
of the operator $\Dc-2x$ in $\Pc_-$ is zero. However
$(\Dc-2x)\Pc_{n-}\subset \Pc_{(n+1)-}.$ Therefore
$(\Dc-2x)\Pc_{n-}= \Pc_{(n+1)-}.$ Finally, 
$$
\la e^{-x^2},\Pc_-\ra_2\ =\
\la e^{-x^2},e^{-x^2}+(\Dc-2x)\Pc_-\ra_2=0
$$ 
and $e^{-x^2}$ belongs to the radical of 
$\la\,.\,,\,.\,\ra_2.$ Since the module $\Pc_-$
is irreducible, the radical has to coincide with the whole
$\Pc_-.$
\qed

Taking real $k,$ the forms  $\la\,.\,,\,.\,\ra_{re}$ 
$\la\,.\,,\,.\,\ra_{im}$ are
positive on $\rit[x]e^{-x^2}$ and
$\rit[\lambda]e^{\lambda^2}.$
The Plancherel formula allows us to complete the function spaces
extending the Fourier transforms $\Fc_{re},\Fc_{im}$
to the spaces of square integrable real-valued functions with
respect to the ``Bessel measure'': 
$$
\Lc^2(\rit,|x|^{2k}dx)\to \Lc^2(i\rit,|\lambda|^{2k}d\lambda)
\to \Lc^2(\rit,|x|^{2k}dx).
$$
The inversion and Plancherel formulas remain valid.

Here we assume that $k>-\frac12.$ Let us discuss the case of
negative half-integers.

\section{Finite-dimensional case}

Let $k=-n-\frac12$ for 
$n\in\zit_+.$ Then  $V_{2n+1}=\Pc/(x^{2n+1})$ is
an irreducible representation of $\HH''$. 
The elements of $V_{2n+1}$ can be identified with polynomials
of degree less than  $2n+1$.

\begin{theorem} 
\label{Clsf}
Finite-dimensional representations of 
$\HH''$ exist only as $k=-n-1/2$ for
$n\in \zit_+.$ 
Given such $k,$ the algebra $\HH''$ has
a unique finite-dimensional irreducible representation
up to isomorphisms, namely, $V_{2n+1}.$
\end{theorem}
{\it Proof.} We will use that 
\begin{align}\label{hcom}
 & [h,x]=x,\ [h,\partial]=-\partial\hbox{\ \ for\ }
h=(x\partial+\partial x)/2.
\end{align}
It readily follows from the defining relations of $\HH''.$
Actually (\ref{hcom}) determines a super
Lie algebra. One may use a general theory of such Lie algebras. 
However in this particular case a reduction to $sl_2$ is
more than sufficient.

Note that $h$
is $x\dx+k+1/2$ in the polynomial
representation. Since the latter is faithful,
(\ref{hcom})  is exactly the claim that  
$x,\Dc$ are homogeneous operators 
of degree $\pm 1,$ which is obvious.
 
We will employ that $e=x^2,$  $f=-\partial^2/4,$ and $h$ 
satisfy the defining relations
of $sl_2(\cit).$ Namely, $[e,f]=h$ because
$$
[\partial^2,x^2]=[\partial^2,x]x+x[\partial^2,x]=
2\partial x+x (2\partial),
$$
and the relations $[h,e]=2e,[h,f]=-2f$ readily result 
from (\ref{hcom}). Cf. Section 6.

Let $V$ be a finite-dimensional representation
of $\HH''.$ 
Then the subspaces $V^0,V^1$ of $V$ 
formed respectively by $s$-invariant and $s$-ani-invariant vectors
are preserved by  $h,$  $e,$ and $f.$ So they are
$sl_2(\cit)$-modules. One gets
$$
\partial x =h+k+1/2,\ x\partial =h-k-1/2 \hbox{\ \ in\ } V^0,
$$
and the other way round in $V^1.$

Let us check that $k\in -1/2- \zit_+.$
All $h$-eigenvalues in $V$ are integers 
thanks to the general theory of finite-dimensional 
representations of $sl_2(\cit)$.
We pick a nonzero $h$-eigenvector $v\in V$ 
with the maximal possible eigenvalue $m.$
Then $m\in \zit_+$ (the theory of $sl_2$) and 
$x(v)=0$ because the latter is an $h$- eigenvector with
the eigenvalue $m+1.$ Hence $\partial x (v)=0,$ 
$m+k+1/2=0,$ and $k=-1/2-m.$ 

Let $U^0$ be a nonzero irreducible 
$sl_2(\cit)-$submodule of $V^0.$ The spectrum of $h$ in $U^0$ is 
$\{\,-n,-n+2,\ldots,n-2,n\,\}$ for an integer $n\ge 0.$
Let $v_l\neq 0$ be an $h$-eigenvector with the eigenvalue 
$l.$ If $e(v)=0$ then $v=cv_{n}$
for a constant $c,$ and if $f(v)=0$ then
$v=cv_{-n}.$  

Let us check that $\partial x (v_{n})=0,$ 
$x\partial (v_{-n})=0,$ and  
$$
\partial x(v_l)\neq 0\hbox{\ for\ } l\neq n, 
\ \  x\partial(v_l)\neq 0 \hbox{\ for\ } l\neq -n.
$$ 

Both operators, $\partial x$ and $x\partial,$ obviously
preserve $U^0:$ 
$$
\partial x(v_l)=(l+k+1/2)v_l,\ x\partial (v_l)=(l-k-1/2)v_l.  
$$
Hence, 
$$
\partial^2 x^2(v_l)=
((\partial x)^2+(1-2k)(\partial x))(v_l)=
(l+k+1/2)(l-k+3/2)v_l.
$$
Setting $l=n,$ we get that $(n+k+1/2)(n-k+3/2)=0$
and $k=-1/2-n,$ because 
$k<0$ and  $n-k+3/2>0.$  Thus $\partial x (v_n)=0.$
The case of $x\partial$ is analogous.

The next claim is that $x(v_n)=0,\ \partial(v_{-n})=0.$ Indeed, 
$x(v')=0$ and $\partial(v')=0$ for  $v'=x(v_n).$ Therefore 
$$
0=[\partial,x](v')=
(1+2ks)(v')=(1-1+n)v'=nv'.
$$
This means that either $v'=0$ or $n=0.$ In the latter case,
$v'$ is proportional to $v_0$ and therefore $v'=x(v_0)=0$
as well. 
Similarly, $\partial(v_{-n})=0.$ 

Now we use the formula 
$$
\partial (x^2(v_l))=x(2+x\partial)(v_l)=
(2+l-k-1/2)x(v_l)=(2+l+n)x(v_l),
$$
and get that $x(v_l)\in \partial(U^0)$ for any $-n\le l\le n.$
Hence $U=U^0+\partial(U^0)$ 
is $x$-invariant. It is obviously $\partial$-invariant
and $s$-invariant ($\Leftarrow$ $\partial(V^0)\subset V^1$).
Also the sum is direct.

Finally, $U$ is a $\HH''$-module and 
has to coincide with $V$ because the latter was assumed to be
irreducible. The above formulas are sufficient to establish
a $\HH''$- isomorphism $U\simeq V_{2n+1}.$ 
Explicitly, the $h$-eigenvectors $x^i(v_{-n})\in U$ 
will be identified 
with the monomials $x^i\in V_{2n+1}.$
\qed

Let us discuss the Hankel transform
and related structures in the case of $V_{2n+1}.$ 
We follow Section 7.

{\bf (a) Form.} To make $\st$ ``inner'' we 
have to construct a non-degenerate bilinear pairing
$(u,v)$ on $V_{2n+1}$ such that $(Hu,v)=(u,H^\st v).$
Here it is:
\begin{equation}
\forall f,g\in V_{2n+1} \mbox{  set  } 
(f,g)=\Res (f(x)g(x)x^{-2n-1}),
\mbox{  where } \Res (\sum a_ix^i) = a_{-1}.
\end{equation}
The pairing is non-degenerate, because if $f=ax^l +$  lower order terms,
where $a\ne 0$ and $0\leq l\leq 2n,$ then $(f,x^{2n-l})=a$. 

We introduce a scalar product $(f,g)_0=\Res(fg)$
for polynomials in terms of $x$ and $x^{-1}.$ 
Denoting 
the conjugate of an operator $A$ with respect to this pairing
by $A^{\circ},$
$$
s^\circ=-s,\ x^\circ=x,\ \dx^{\circ}=-\dx, \hbox{\ \, and\ \,}
x^{2n+1}\Dc^{\circ}x^{-2n-1}=-\Dc.
$$
The relation $x^\circ=x$ is obvious. Concerning
$\dx^{\circ}=-\dx,$ it 
follows from the property $\Res (df/dx)$=0 for any 
polynomial $f(x,x^{-1}).$ The formula $s^{\circ}=-s$
results from $\Res(s(f(x))=-\Res(f(x)).$  

Switching from ${}^\circ$ to ${}^\st,$
we have 
$$ (sf,g)=\Res(s(f)gx^{-2n-1})=-\Res(fs(g)s(x^{-2n-1}))=
-(-1)^{2n+1}(fs(g)x^{-2n-1})=(f,s(g)).
$$
Finally,
\begin{align}
\Dc^{\circ}&=(\dx+\frac{k}{x}(1-s))^{\circ}=-\dx+(1+s)\frac{k}{x},\notag\\
x^{-2k}\Dc^{\circ}x^{2k}&= -\dx -\frac{2k}{x}+\frac{k}{x}(1+s)=
-\dx + \frac{k}{x}(-1+s)= -\Dc.
\end{align}
The first equality on the second line holds because $\frac{k}{x}x^{2k}=
kx^{-2n-2}$ is an even function, and thus it commutes with the action of
$s$. Finally
$$ (\Dc f,g)=(\Dc f,x^{2k}g)_0=
(f,x^{2k}x^{-2k}\Dc^{\circ}x^{2k}(g))_0
=-(f,x^{2k}\Dc g)_0=-(f,\Dc g).\quad\square$$

{\bf (b) Gaussian.}  The Gaussian does not exist in polynomials
but of course can be introduced as a power series
$e^{x^2}=\sum_{m=0}^{\infty}(x^2)^m/m!$ in the algebra of 
formal series $\cit[[x]],$ a completion
of the polynomial representation.
The conjugation by this series induces $\tau$ on $\HH''.$  
Its inverse is $e^{-x^2}=\sum_{m=0}^{\infty}(-x^2)^m/m!.$
The multiplication by the Gaussian
does not preserve the space of polynomials but is well-defined on
$V_{2n+1}$ because $\forall f\in V_{2n+1}$ we have $x^mf=0$ 
for $m\geq 2n+1.$ 
Finally, 
$$\gamma^{\pm}=\sum_{m=0}^{2n}(\pm x^2)^m/m!.
$$

{\bf (c) Hankel transform.}
The operator $\Dc$ is nilpotent in $V_{2n+1}$ because it lowers the
degree of $f\in V_{2n+1}$ by one. Therefore the operators 
$\exp (\pm \Dc^2/4)\in \cit[[\Dc]]$ are well-defined
in this representation as well as the Gaussians.
It suffices to take 
$\sum_{m=0}^{2n}(\pm (D/2)^2)^m/m!$.
Thus we may set 
\begin{equation}\label{trunc-H}
  \Fb = e^{x^2}e^{\frac{\Dc^2}{4}}e^{x^2}
\hbox{ \ in\ } V_{2n+1}.
\end{equation}

\begin{proposition}\label{Htran}
The map $\Fc$ is a Hankel transform on $V_{2n+1}$, i.e. 
$\Fb\circ\Dc=-\Fb\circ 2x$, $\Fb\circ 2x=\Dc\circ\Fb$, 
$\Fb s=\Fb s$. 
These relations fix it uniquely up to proportionality.
\end{proposition}
{\it Proof.} 
We already know that $\Fb$ is a Hankel transform
(the previous section).
If $\widetilde{\Fb}$ is another one
then the ratio $\widetilde{\Fb}\Fb^{-1}$ commutes
with  $x,\Dc$, and $s$ because of the very definition.
Since $V_{2n+1}$ is irreducible (and finite dimensional)
we get that $\widetilde{\Fb}$ is proportional to $\Fb.$
\qed

\section{Truncated Bessel functions}

Recall that $\Dc$ has only one eigenvalue in $V_{2n+1}$,
namely, $0.$ Therefore we cannot define the $\psi_\lambda$
as an eigenfunction
of $\Dc$ in $V_{2n+1}$ any longer. 
Instead, it will be introduced as the 
kernel of the Hankel transform. 

Any linear operator 
$A:V_{2n+1}^x\rightarrow V_{2n+1}^{\lambda}$ (the
upper index indicates the  variable) is a matrix. It means
that
\begin{align}
&A(f)(\lambda)=(f(x),\alpha(x,\lambda))=
\Res(f(x)\alpha(x,\lambda)x^{-2n-1}),\hbox{\ where \ }\notag\\
&\alpha(x,\lambda)=\sum_{l,m=0}^{2n}c_{l,m}x^l\lambda^m
=\sum_{l=0}^{2n}x^{2n-l}A(x^l).
\end{align}
So here the kernel $\alpha(x,\alpha)$ is uniquely defined by $A$
and vice versa.

The {\it truncated $\psi$-function} is the kernel of $\Fb:$  
\begin{equation}
\Fb(f)(\lambda)=(f(x),\psib_\lambda(x))=
\Res(f(x)\psib_\lambda(x)x^{-2n-1}).
\end{equation}

There is a somewhat different approach.
Let us use that 
the relations from Lemma~\ref{Htran} determine $\Fb$
uniquely up to
proportionality. These relations are equivalent to the
following properties of $\psib_\lambda(x):$
\begin{equation}\label{tef}
\Dc\psib_\lambda(x)=2\lambda\psib_\lambda(x) 
\mod (x^{2n+1},\lambda^{2n+1}),\
\psib_\lambda(x)=\psib_x(\lambda),\ 
\psib_\lambda(s(x))=\psib_{s(\lambda)}(x).
\end{equation}

Let us solve the first equation. Setting
$\psil(x)=\sum_{l,m=0}^{2n}c_{l,m}x^l\lambda^m,$ 
\begin{align} 
&\sum_{l=1,m=0}^{l=2n,m=2n} c_{l,m}(l+(1-(-1)^l)
(-n-\frac12))x^{l-1}\lambda^m\notag\\
=&\sum_{l=0,m=0}^{l=2n,m=2n-1} 2c_{l,m}x^l\lambda^{m+1} \mod
(x^{2n+1},\lambda^{2n+1}),\notag\\
&c_{l,m}=\frac{2}{l+(1-(-1)^l)(-n-1/2)}c_{l-1,m-1}
\mbox{ for } 2n\ge l>0,\ 2n\ge m>0,\notag \\
&\mbox{where } c_{l,0}=0=c_{2n,m} \mbox{ for } l>0,\ m<2n.
\label{fceq}
\end{align}

Using the $x\leftrightarrow \lambda$ symmetry, we conclude
that $c_{l,0}=0=c_{0,l}$ for nonzero $l$ and
$c_{l,m}=0$ for $l\neq m.$ Thus 
$$
\psib_\lambda(x)=g_n(\lambda x)\mbox{ for }
g_n=\sum_{l=0}^{2n}c_{l}t^l,\ c_l=c_{l,l},
$$
where the coefficients are given by (\ref{fceq}).

Finally, $g_n(t)=f_n(t)+(1/2)df_n/dt$ for the 
{\it truncated Bessel function}
$f_n(t)=\sum_{m=0}^{n}c_{2m}t^{2m}$ which is an even
solution of the truncated Bessel equation (cf. Section 1):
\begin{equation}\label{btrun}
\frac{d^2f}{dt^2}(t) +
2k\frac{1}{t}\frac{df}{dt}(t)-4f(t)=0\, \mod\, (t^{2n}),\ 
k=-n-1/2.
\end{equation}
This equation is sufficient to determine the coefficients
of $f_n$ uniquely for any constant term $c_{0}=c_{0,0}.$ 
They are given by the same formula (\ref{feqn}) till
$c_{2n}$ up to proportionality. This can be checked directly
using explicit formulas which will be discussed next.

Still $c_0$ remains arbitrary.
Recall that $\psib$ was initially introduced
as the kernel of $\Fb.$ So it comes with its own normalization.
Let us calculate  its $c_{0}.$
One gets:
\begin{align}
&\Fb(e^{-x^2})=e^{\lambda^2}\exp(\Dc^2/4)(e^{\lambda^2}e^{-\lambda^2})
=e^{\lambda^2}\exp(\Dc^2/4)(1)=e^{\lambda^2},\mbox { so }\notag\\
&\Fb(1-\frac{x^2}{1!}+\frac{x^4}{2!}+\cdots
+(-1)^n\frac{x^{2n}}{n!})=1+\frac{\lambda^2}{1!}+\cdots
\frac{\lambda^{2n}}{n!}\label{ftrc}.
\end{align}
Here the transform of $1$ is proportional
to $\lambda^{2n}$ since the latter
has to be an eigenfunction of $\lambda,$ i.e. the solution
of the equation $\lambda \Fb(1)=0$ in $V_{2n+1}^\lambda.$ 
Similarly, $\Fb(x^{l})=(\Dc^\lambda/2)^{l}\Fb(1)$
is proportional to $\lambda^{2n-l}$ for $0\le l\le 2n.$
Thus (\ref{ftrc}) leads to the relations  
\begin{align}\label{treven}
&\Fb(x^{2m})=
(-1)^m\frac{m!}{(n-m)!}\lambda^{2n-2m}.
\end{align}
For instance, $\Fb(x^{2n})=
(-1)^n n!.$ This is exactly the coefficient $c_0$ above.

We obtain that the normalization serving the truncated
Hankel transform is
\begin{equation}\label{fnorm}
\psib_\lambda(0)=-n!,\ c_0=g_n(0)=f_n(0)=-n!.
\end{equation}

Formula (\ref{treven}) also results in 
\begin{align}\label{trodd}
&\Fb(x^{2m+1})=\Fb(x(x^{2m}))=(\Dc/2)\Fb(x^{2m})=
(-1)^m\frac{m!}{(n-m-1)!}\lambda^{2n-2m-1}.
\end{align}
Substituting,
\begin{align}
&\psib_\lambda(x)=\sum_{m=0}^{n}x^{2n-2m}\Fb(x^{2m})+
\sum_{m=0}^{n-1}x^{2n-2m-1}\Fb(x^{2m+1})
=f_n(x\lambda)+\frac12 f'_n(x\lambda)\hbox{\ for\ } \notag\\
&f_n(t)=\sum_{m=0}^{n}\frac{(-1)^mm!}{(n-m)!}t^{2n-2m}=
\sum_{m=0}^{n}\frac{(-1)^{n-m}(n-m)!}{m!}t^{2m}.\label{finf}
\end{align}
It is exactly the solution of (\ref{btrun}) with the 
{\it truncated normalization} $f_n(0)=(-1)^n n!.$

{\bf Truncated inversion.}
Concluding the consideration
of the case $k=-n-\frac12$ for  $n\in\zit_+,$ 
let us discuss the inversion. We have the 
following transformations and scalar products:
\begin{align}
\Fb_+:\cit[x]/(x^{2n+1})\rightarrow 
\cit[\lambda]/(\lambda^{2n+1}),&
\qquad\Fb_+(f)=\Res(f(x)\psib_\lambda(x)x^{-2n-1}),\notag\\
\Fb_-:\cit[\lambda]/(\lambda^{2n+1})
\rightarrow \cit[x]/(x^{2n+1}),&
\qquad\Fb_-(f)=
\Res(f(\lambda)\psib_x(-\lambda)\lambda^{-2n-1}).\notag\\
\la f,g\ra_+=\Res(f(x)g(x)x^{-2n-1}),&\qquad f,g\in \cit[x]/(x^{2n+1});
\notag\\
\la f,g\ra_-=\Res(f(-\lambda)g(\lambda)\lambda^{-2n-1}),&\qquad f,g\in
\cit[\lambda]/(\lambda^{2n+1}).
\end{align} 

Here $\Fb_+(f)=\Fb(f)=\widehat{f}$ in the notation above.
The transform $\Fb_-(f)$ coincides with $\Fb^\lambda_+(f)$ for even
$f(\lambda)$ and with $-\Fb^\lambda_+(f)$ for odd $f(\lambda).$

We can follow the ``analytic'' case and check that 
$\Fb_-\circ\Fb_+$ commutes with $\Dc, x, s.$ Hence it
is the multiplication by a constant thanks to the irreducibility
of $V_{2n+1}.$ The constant is $\Fb_-\circ\Fb_+(1)$ and can be
readily calculated. It is equally simple to calculate all 
$\Fb_-\circ\Fb_+(x^l)$ using (\ref{treven}) and (\ref{trodd}).
For instance,
\begin{align}
&\Fb_-\circ\Fb_+(x^{2m})=\Fb_-\left(\frac{(-1)^m m!}{(n-m)!}
\lambda^{2n-2m}\right)=\notag\\
&=\frac{(-1)^mm!}{(n-m)!}
\frac{(-1)^{n-m}(n-m)!}{m!}x^{2m}=(-1)^nx^{2m}.\notag
\end{align}  
Thus the truncated inversion reads: 
$$
\Fb_-\circ\Fb_+\ =\ (-1)^n\,\hbox{id}\ =\ \Fb_+\circ\Fb_-.
$$

Concerning the Plancherel formula, 
we may use the proportionality of the forms
$\la f,g\ra_+$ and $\la \widehat{f},\widehat{g}\ra_-$
for $f,g\in V_{2n+1}$ and their transforms 
$\widehat{f}=\Fb(f), \widehat{g}=\Fb(g).$ It results
from the irreducibility of $ V_{2n+1}.$ A direct calculation
is simple as well. Let
\begin{align}
&\la f,f\ra_+=\la f,f\ra= \sum_{l=0}^{2n} a_l a_{2n-l} 
\hbox{ for } f=\sum_{l=0}^{2n} a_l x^l,\notag\\
&\la g,g\ra_-= \sum_{l=0}^{2n} (-1)^l b_l b_{2n-l} 
\hbox{ for } g=\Fb(f)=\sum_{l=0}^{2n} b_l \lambda^l.
\end{align}
It is easy to check that 
$$
b_lb_{2n-l}=(-1)^{l+n} a_l a_{2n-1}.
$$
Indeed, using (\ref{treven}) and (\ref{trodd}): 
\begin{align}
&b_{2m} b_{2n-2m}=
(-1)^m a_{2m}
\frac{m!}{(n-m)!}(-1)^{n-m}a_{2n-2m}\frac{(n-m)!}{m!}\notag\\
&=(-1)^n a_{2m}a_{2n-2m},\notag\\
&b_{2m+1} b_{2n-2m-1}=
(-1)^m a_{2m+1}
\frac{m!}{(n-m-1)!}(-1)^{n-m-1}a_{2n-2m-1}\frac{(n-m-1)!}{m!}
\notag\\
&=(-1)^{n-1} a_{2m+1} a_{2n-2m-1}.\notag
\end{align}
We get the truncated Plancherel formula:
$$ 
\la \widehat{f},\widehat{g}\ra_-\ =\ (-1)^n\,\la f,g\ra_+.
$$
The above consideration proves the coincidence for $f=g,$
i.e. for the corresponding quadratic forms. It is of course sufficient.


\begin{thebibliography}{MTV}

\bibitem[C1]{C1}
Cherednik, I.:
Difference Macdonald-Mehta conjectures.
IMRN {\bf 10}, 449--467 (1997).
\bibitem[C2]{C2}
Cherednik, I.:
One-dimensional double Hecke algebras and Gaussians,
CIME (2000).
\bibitem[D]{D} 
Dunkl, C.F.: Differential-difference operators associated to
reflection groups, Trans. AMS. {\bf 311}, 167--183 (1989).
\bibitem[J]{J}
Jeu, M.F.E. de:
The Dunkl transform, Invent. Math. {\bf 113},  147--162 (1993).
\bibitem[L] {L}
Luke, J.: Integrals of Bessel functions, McGraw-Hill Book Company,
New York-Toronto-London (1962).
\bibitem[O] {O}
Opdam, E.M.: Dunkl operators, Bessel functions and the 
discriminant of a finite Coxeter group,
Comp. Math. {\bf 85}, 333--373 (1993).

 
\end{thebibliography}
\end{document}